\documentclass{llncs}
\usepackage{amssymb}
\usepackage{pstricks}
\usepackage{pst-node}
\usepackage{color}
\usepackage{epsf}
\usepackage{latexsym}

\begin{document}

\title{Toward a Unit Distance Embedding for the Heawood graph}

\author{
Mitchell A. Harris\inst{1}}

\institute{
Harvard Medical School\\
Department of Radiology\\
Boston, MA 02114, USA\\
\email{harris.mitchell@mgh.harvard.edu}
}

\maketitle

\begin{abstract}

The unit distance embeddability of a graph, like planarity, involves a
mix of constraints that are combinatorial and geometric. We construct
a unit distance embedding for $H-e$ in the hope that it will lead to
an embedding for $H$. We then investigate analytical methods for a
general decision procedure for testing unit distance embeddability.

\end{abstract}

\section{Introduction}

Unit distance embedding of a graph is an assignment of coordinates in
the plane to the vertices of a graph such that if there is an edge
between two vertices, then their respective coordinates are exactly
distance 1 apart. To bar trivial embeddings, such as for bipartite
graphs having all nodes in one part located at $(0,0)$, for the other
part at $(0,1)$, it is also required that the embedded points be
distinct. There is no restriction on edge crossing. A graph is said to
be unit distance embeddable if there exists such an embedding (with
the obvious abbreviations employed, such as UD embeddable, a UDG,
etc).

For some well-known examples, $K_4$ is not UDG but $K_4 - e$ is. The
Moser spindle is UDG, and is also 4 colorable \cite{MM61}, giving the
largest currently known lower bound to the Erd\"os colorability of the
plane \cite{BE51}. The graph $K_{2,3}$ is not because from two given
points in the plane there are exactly two points of distance 1, but
the graph wants three; removing any edge allows a UD embedding.

The property UD is hereditary so by the Graph Minor theorem, the
property has a finite number of forbidden minors, and by the previous
paragraph, $K_4$ and $K_{2,3}$ are two of them.

The Petersen graph $P$ is a UDG, as well as $K_2^n$.

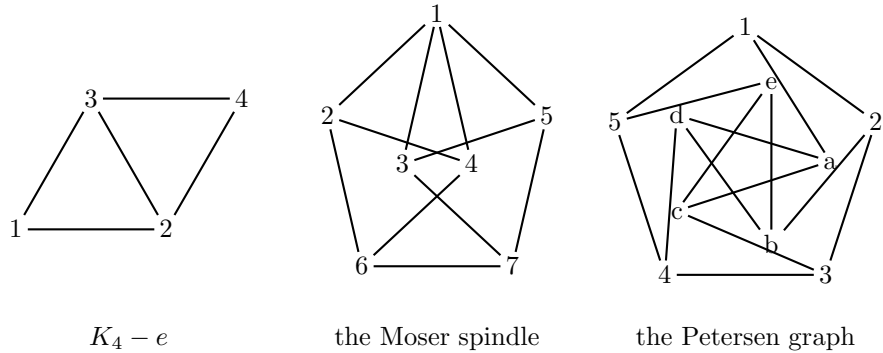
\psset{unit=2cm, nodesep=2pt}

\begin{center}
$
\begin{array}{ccc}
\begin{pspicture}(-1,-1)(1,1)
\rput(-0.75, -0.435){\Rnode{v1}{1}}
\rput(0.25, -0.435){\Rnode{v2}{2}}
\rput(-0.25, 0.435){\Rnode{v3}{3}}
\rput(0.75, 0.435){\Rnode{v4}{4}}
\ncline{v1}{v2}\ncline{v1}{v3}\ncline{v2}{v3}\ncline{v2}{v4}\ncline{v3}{v4}
\end{pspicture}
&
\begin{pspicture}(-1,-1)(1,1)
\rput(0., 1.0){\Rnode{v1}{1}}
\rput(-0.728714, 0.32){\Rnode{v2}{2}}
\rput(-0.228714, 0.0){\Rnode{v3}{3}}
\rput(0.228714, 0.0){\Rnode{v4}{4}}
\rput(0.728714, 0.32){\Rnode{v5}{5}}
\rput(-0.5, -0.68){\Rnode{v6}{6}}
\rput(0.5, -0.68){\Rnode{v7}{7}}
\ncline{v1}{v2}\ncline{v1}{v3}\ncline{v1}{v4}\ncline{v1}{v5}
\ncline{v2}{v4}\ncline{v3}{v5}
\ncline{v2}{v6}\ncline{v4}{v6}
\ncline{v3}{v7}\ncline{v5}{v7}
\ncline{v6}{v7}
\end{pspicture}
&
\begin{pspicture}(-1,-1)(1,1)
\rput(0.0, 0.911){\Rnode{v1}{1}}
\rput(0.866, 0.282){\Rnode{v2}{2}}
\rput(0.534, -0.737){\Rnode{v3}{3}}
\rput(-0.534, -0.737){\Rnode{v4}{4}}
\rput(-0.866, 0.282){\Rnode{v5}{5}}
\rput(0.563, 0.0){\Rnode{va}{a}}
\rput(0.174, -0.536){\Rnode{vb}{b}}
\rput(-0.455, -0.331){\Rnode{vc}{c}}
\rput(-0.455, 0.331){\Rnode{vd}{d}}
\rput(0.174, 0.536){\Rnode{ve}{e}}
\ncline{v1}{v2}\ncline{v2}{v3}\ncline{v3}{v4}\ncline{v4}{v5}\ncline{v5}{v1}
\ncline{va}{vc}\ncline{vb}{vd}\ncline{vc}{ve}\ncline{vd}{va}\ncline{ve}{vb}
\ncline{v1}{va}\ncline{v2}{vb}\ncline{v3}{vc}\ncline{v4}{vd}\ncline{v5}{ve}
\end{pspicture}
\\
K_4 - e & \hbox{the Moser spindle} & \hbox{the Petersen graph}
\end{array}
$
\end{center}

A UDG is {\it rigid} if there are only a finite number of UD
embeddings, that is, one embedding cannot be transformed continuously
to any other.

Now consider the Heawood graph $H$, also known as the point-line graph
of the Fano plane, the (3,6) cage, or the smallest cubic graph of
girth 6, and specified by LCF notation as $(5, -5)^7$ or the
difference set $\{1,2,4\} \bmod 14$. It has the following non-UD
embedding:

\begin{center}
$
\begin{array}{c}
\begin{pspicture}(-1,-1)(1,1)
\rput(0.900969, 0.433884){\Rnode{vb}{b}}
\rput(0.62349, 0.781831){\Rnode{v2}{2}}
\rput(0.222521, 0.974928){\Rnode{va}{a}}
\rput(-0.222521, 0.974928){\Rnode{v1}{1}}
\rput(-0.62349, 0.781831){\Rnode{vg}{g}}
\rput(-0.900969, 0.433884){\Rnode{v7}{7}}
\rput(-1., 0){\Rnode{vf}{f}}
\rput(-0.900969, -0.433884){\Rnode{v6}{6}}
\rput(-0.62349, -0.781831){\Rnode{ve}{e}}
\rput(-0.222521, -0.974928){\Rnode{v5}{5}}
\rput(0.222521, -0.974928){\Rnode{vd}{d}}
\rput(0.62349, -0.781831){\Rnode{v4}{4}}
\rput(0.900969, -0.433884){\Rnode{vc}{c}}
\rput(1., 0){\Rnode{v3}{3}}
\ncline{v1}{va}\ncline{va}{v2}\ncline{v2}{vb}\ncline{vb}{v3}\ncline{v3}{vc}\ncline{vc}{v4}\ncline{v4}{vd}
\ncline{vd}{v5}\ncline{v5}{ve}\ncline{ve}{v6}\ncline{v6}{vf}\ncline{vf}{v7}\ncline{v7}{vg}\ncline{vg}{v1}
\ncline{v1}{vc}\ncline{v2}{vd}\ncline{v3}{ve}\ncline{v4}{vf}\ncline{v5}{vg}\ncline{v6}{va}\ncline{v7}{vb}
\end{pspicture}
\\\\
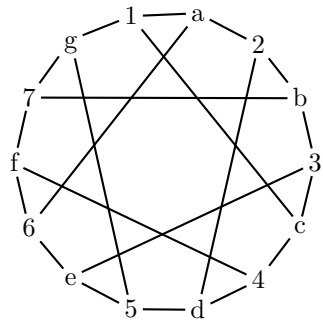
\hbox{The Heawood graph $H$}  
\end{array}
$
\end{center}

\section{The Construction}
Consider first the Heawood graph with two adjacent vertices
removed. We will show that this graph is UDE and then add back in the
two vertices (but not the mutual incident edge). The graph $H -
\{1,a\}$ ((a) in the figure) is isomorphic to the M\"obius ladder
$M_4$ with a vertex inserted on each 'rung' (b):

\begin{center}
$  
\begin{array}{ccc}
\begin{pspicture}(-1,-1)(1,1)
\rput(0.900969, 0.433884){\Rnode{vb}{\red{b}}}
\rput(0.62349, 0.781831){\Rnode{v2}{2}}
\rput(-0.62349, 0.781831){\Rnode{vg}{g}}
\rput(-0.900969, 0.433884){\Rnode{v7}{\red{7}}}
\rput(-1., 0){\Rnode{vf}{\red{f}}}
\rput(-0.900969, -0.433884){\Rnode{v6}{6}}
\rput(-0.62349, -0.781831){\Rnode{ve}{\red{e}}}
\rput(-0.222521, -0.974928){\Rnode{v5}{\red{5}}}
\rput(0.222521, -0.974928){\Rnode{vd}{\red{d}}}
\rput(0.62349, -0.781831){\Rnode{v4}{\red{4}}}
\rput(0.900969, -0.433884){\Rnode{vc}{c}}
\rput(1., 0){\Rnode{v3}{\red{3}}}
\ncline{v2}{vb}\ncline{vb}{v3}\ncline{v3}{vc}\ncline{vc}{v4}\ncline{v4}{vd}
\ncline{vd}{v5}\ncline{v5}{ve}\ncline{ve}{v6}\ncline{v6}{vf}\ncline{vf}{v7}\ncline{v7}{vg}
\ncline{v2}{vd}\ncline{v3}{ve}\ncline{v4}{vf}\ncline{v5}{vg}\ncline{v7}{vb}
\end{pspicture}
&
\qquad\phantom{x}\qquad
&
\begin{pspicture}(-1,-1)(1,1)
\rput(0, -1){\Rnode{v5}{\red{5}}}
\rput(-0.707, -0.707){\Rnode{vd}{\red{d}}}
\rput(-1., 0){\Rnode{v4}{\red{4}}}
\rput(-0.707, 0.707){\Rnode{vf}{\red{f}}}
\rput(0, 1){\Rnode{v7}{\red{7}}}
\rput(0.707, 0.707){\Rnode{vb}{\red{b}}}
\rput(1., 0){\Rnode{v3}{\red{3}}}
\rput(0.707, -0.707){\Rnode{ve}{\red{e}}}
\rput(0.2, -0.2){\Rnode{vg}{g}}
\rput(-0.2, -0.2){\Rnode{vc}{c}}
\rput(-0.2, 0.2){\Rnode{v2}{2}}
\rput(0.2, 0.2){\Rnode{v6}{6}}
\ncline{v2}{vb}\ncline{vb}{v3}\ncline{v3}{vc}\ncline{vc}{v4}\ncline{v4}{vd}
\ncline{vd}{v5}\ncline{v5}{ve}\ncline{ve}{v6}\ncline{v6}{vf}\ncline{vf}{v7}\ncline{v7}{vg}
\ncline{v2}{vd}\ncline{v3}{ve}\ncline{v4}{vf}\ncline{v5}{vg}\ncline{v7}{vb}
\end{pspicture}
\\\\
\hbox{a) $H - \{1,a\}$} &
&
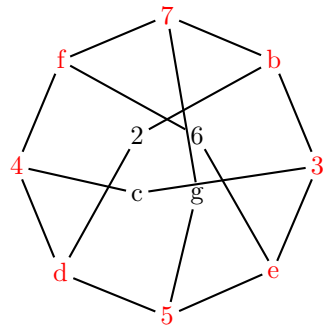
\hbox{b) as a ladder with interposed vertices on the rungs}
\end{array}
$
\end{center}

The difficulty with UD embedding this modified ladder is that opposing
vertices are mutually more than distance 2 apart. So we transform a
smaller graph and then build back up. First, it is easy to UD embed
the nine points that make up the 8-cycle and one rung on a
square. 'Folding' over the rung puts the opposite corners closer
together, and perturbing $f$, $d$, and $4$ a little preserves
distinctness (a). The rest of the rungs can then be added since the
end vertices are now all within distance 2 (b). The last two points of
$H$ can now be added in (c) giving a UDE of $H - $.

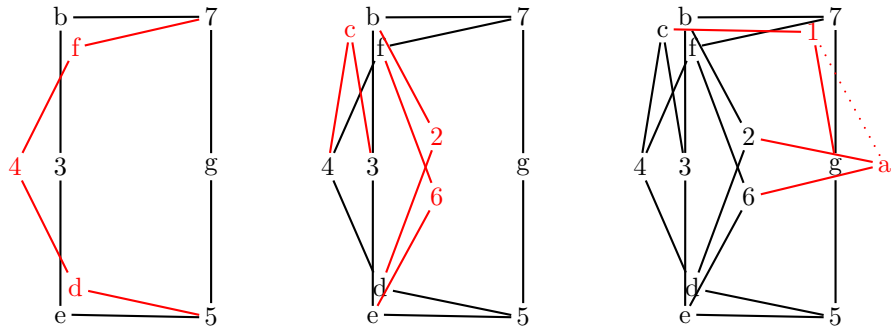
\begin{center}
$
\begin{array}{ccc}
\begin{pspicture}(-1,-1)(1,1)
\rput(0.5, -1){\Rnode{v5}{5}}
\rput(-0.4, -0.8){\Rnode{vd}{\red{d}}}
\rput(-0.8, 0){\Rnode{v4}{\red{4}}}
\rput(-0.4, 0.8){\Rnode{vf}{\red{f}}}
\rput(0.5, 1){\Rnode{v7}{7}}
\rput(-0.5, 1){\Rnode{vb}{b}}
\rput(-0.5, 0){\Rnode{v3}{3}}
\rput(-0.5, -1){\Rnode{ve}{e}}
\rput(0.5, 0){\Rnode{vg}{g}}
\ncline{v7}{vb}\ncline{vb}{v3}\ncline{v3}{ve}\ncline{ve}{v5}
\ncline[linecolor=red]{v5}{vd}\ncline[linecolor=red]{vd}{v4}\ncline[linecolor=red]{v4}{vf}\ncline[linecolor=red]{vf}{v7}
\ncline{v7}{vg}\ncline{vg}{v5}
\end{pspicture}
&
\begin{pspicture}(-1,-1)(1,1)
\rput(0.5, -1){\Rnode{v5}{5}}
\rput(-0.45, -0.8){\Rnode{vd}{d}}
\rput(-0.8, 0){\Rnode{v4}{4}}
\rput(-0.45, 0.8){\Rnode{vf}{f}}
\rput(0.5, 1){\Rnode{v7}{7}}
\rput(-0.5, 1){\Rnode{vb}{b}}
\rput(-0.5, 0){\Rnode{v3}{3}}
\rput(-0.5, -1){\Rnode{ve}{e}}
\rput(0.5, 0){\Rnode{vg}{g}}
\rput(-0.65, 0.9){\Rnode{vc}{\red{c}}}
\rput(-0.075, 0.2){\Rnode{v2}{\red{2}}}
\rput(-0.075, -0.2){\Rnode{v6}{\red{6}}}
\ncline{v7}{vb}\ncline{vb}{v3}\ncline{v3}{ve}\ncline{ve}{v5}
\ncline{v5}{vd}\ncline{vd}{v4}\ncline{v4}{vf}\ncline{vf}{v7}
\ncline{v7}{vg}\ncline{vg}{v5}
\ncline[linecolor=red]{vb}{v2}\ncline[linecolor=red]{v2}{vd}
\ncline[linecolor=red]{v3}{vc}\ncline[linecolor=red]{vc}{v4}
\ncline[linecolor=red]{ve}{v6}\ncline[linecolor=red]{v6}{vf}
\end{pspicture}
&
\begin{pspicture}(-1,-1)(1,1)
\rput(0.5, -1){\Rnode{v5}{5}}
\rput(-0.45, -0.8){\Rnode{vd}{d}}
\rput(-0.8, 0){\Rnode{v4}{4}}
\rput(-0.45, 0.8){\Rnode{vf}{f}}
\rput(0.5, 1){\Rnode{v7}{7}}
\rput(-0.5, 1){\Rnode{vb}{b}}
\rput(-0.5, 0){\Rnode{v3}{3}}
\rput(-0.5, -1){\Rnode{ve}{e}}
\rput(0.5, 0){\Rnode{vg}{g}}
\rput(-0.65, 0.9){\Rnode{vc}{c}}
\rput(-0.075, 0.2){\Rnode{v2}{2}}
\rput(-0.075, -0.2){\Rnode{v6}{6}}
\rput(0.35, 0.9){\Rnode{v1}{\red{1}}}
\rput(0.825, 0){\Rnode{va}{\red{a}}}
\ncline{v7}{vb}\ncline{vb}{v3}\ncline{v3}{ve}\ncline{ve}{v5}
\ncline{v5}{vd}\ncline{vd}{v4}\ncline{v4}{vf}\ncline{vf}{v7}
\ncline{v7}{vg}\ncline{vg}{v5}
\ncline{vb}{v2}\ncline{v2}{vd}
\ncline{v3}{vc}\ncline{vc}{v4}
\ncline{ve}{v6}\ncline{v6}{vf}
\ncline[linecolor=red]{v6}{va}\ncline[linecolor=red]{va}{v2}
\ncline[linecolor=red]{vg}{v1}\ncline[linecolor=red]{v1}{vc}
\ncline[linecolor=red,linestyle=dotted]{v1}{va}
\end{pspicture}
\\\\
\hbox{The cycle and one rung,}&
\hbox{plus the other three rungs,}&
\hbox{plus the last two points}
\end{array}
$
\end{center}

Though this is essentially a proof by picture, necessitating all the
usual Euclidean caveats about inferring from idiosyncrasies of the
specific diagram, it still follow. Each of the transformations can be
seen to preserve or enforce unit distance and preserve
distinctness. When a vertex and two edges are added, the UD condition
can be satisfied in exactly two ways. The only real choice made here
is for vertices 2, 6, $a$ and 1, and we make those choices (exercise
for the reader) such that things work.

In order to get a UDE of the full graph $H$, the last item to take
care of is the edge between 1 and $a$. Given that the suggested
embedding is highly constrained by the strict placement of the six
initial vertices, and the freedoms in the perturbations and binary
choices for the rest, can things be modified slightly enough so that 1
and $a$ are a unit apart {\em and} vertex embeddings are kept
distinct?

If one could show that there is a configuration where $\overline{1a} <
1$, a configuration another where $\overline{1a} > 1$, {\it and} a
continuous transformation between the two, then we'd have a proof of
the existence of a unit embedding. This is not exactly a constructive
embedding but a proof nonetheless, from which a numerically accurate
embedding can be approximated.



\section{Analytic and Automatic Solutions}
For some graphs there are obvious 'by-hand' proofs or disproofs of
embeddability or the lack thereof. But we also seek a general
algorithm to determine UDE.

An unit distance embedding graph can be modeled by a set of
multinomial equations that express the coordinates of the vertices of
an edge in a distance constraint. For example, if $x_a, y_a$ and $x_b,
y_b$ are the coordinates of an edge between $a$ and $b$, then by the
Euclidean distance formula:
\[
(x_a - x_b)^2 + (y_a - y_b)^2  = 1,
\]
and each edge of a graph produces another such equality constraint.

For a set of non-linear multinomial equations, there is a decision
procedure that, though it doesn't necessarily produce closed-form
coordinates, it will give a yes-no answer to whether the set of
constraints has a solution. Gr\"obner basis completion takes a set of
multinomials and reduces it to a 'minimal' set, such that the minimal
set has the multinomial '1' as its sole member if and only if there is
no solution. If this minimal set, called the Gr\"obner basis, is not
'1', then it is a set of multinomials, from which one attempt to
numerically extract coordinates (variations on multivariate
Newton-Raphson with all their attendant problems of convergence), or
using other methods, attempt to symbolically extract coordinates
(polynomial factoring, root extraction, etc).

For example, $K_4 - e$, with edges $(1,2)$, $(1,3)$, $(2,3)$, $(2,4)$,
$(3,4)$, has the system:
\begin{eqnarray*}
(x_1 - x_2)^2 + (y_1 - y_2)^2  &=& 1\\
(x_1 - x_3)^2 + (y_1 - y_3)^2  &=& 1\\
(x_2 - x_3)^2 + (y_2 - y_3)^2  &=& 1\\
(x_2 - x_4)^2 + (y_2 - y_4)^2  &=& 1\\
(x_3 - x_4)^2 + (y_3 - y_4)^2  &=& 1
\end{eqnarray*}

To oversimplify, the completion algorithm will do a generalization of
Gaussian elimination, eliminating largest common terms between two
(expanded) multinomial equations. For the above system, setting
$(x_1,y_1)$ to $(0,0)$ and $(x_2,y_2)$ to $(1,0)$ to reduce some
processing, we get the following reduced system:

\begin{eqnarray*}
4 y_4^3 &=& 3 y_4 \\
x4 &=& 2 y4^2\\
y3 y4 &=& y4^2\\
4 y3^2 &=& 3\\
2 x3 &=& 1\\
y2 &=& 0\\
x2 &=& 1\\
y1 &=& 0\\ 
x1 &=& 0
\end{eqnarray*}
where, in this instance, it is easy to extract the values of the
coordinates by back-substitution to get:
\begin{eqnarray*}
  (x_1,y_1) &=& (0,0)\\
  (x_2,y_2) &=& (1,0)\\
  (x_3,y_3) &=& ({1\over2}, \pm {\sqrt{3} \over 2})\\
  (x_4,y_4) &=& (0,0) \hbox{ or } ({3 \over 2}, \pm {\sqrt{3} \over 2}) \\
\end{eqnarray*}

This answer has to be checked for duplicate vertex embeddings, but
when $(x_4,y_4) = (0,0)$ is removed, there is still a legal embedding
left.

For the example of $K_{2,3}$, where there are many (continuous)
solutions to the system. Computing by mechanically using a symbolic
algebra package, we get:
\begin{eqnarray*}
  x_5^2 + y_5^2 &=& 1\\
  x_4^2 + y_4^2 &=& 1\\
  x_2 x_5 + y_2 y_5 &=& x_2\\
  x_2 x_4 + y_2 y_4 &=& x_2\\
  x_2^2 + y_2^2 &=& 2 x_2\\
  y_2^2 y_5 &=& y_2 y_5^2\\
  y_2 y_5 + x_5 y_2 y_5 &=& x_2 y_5^2\\
  x_2 y_2 y_5 &=& x_2 y_5^2\\
  y_2^2 y_4 &=& y_2 y_4^2\\
  x_5 y_2 y_4 + y_2 y_5 &=& x_4 y_2 y_5 + y_2 y_4\\
  y_2 y_4 + x_4 y_2 y_4 &=& x_2 y_4^2\\
  x_2 y_2 y_4 &=& x_2 y_4^2\\
  x_5 y_2^2 + 2 y_2 y_5 &=& x_2 y_5^2 + y_2^2\\
  x_4 y_2^2 + 2 y_2 y_4 &=& x_2 y_4^2 + y_2^2\\
  2 y_2 y_4 y_5 + x_4 y_2 y_5^2 &=& y_2 y_5^2 + x_2 y_4 y_5^2\\
  y_2 y_4^2 y_5 &=& y_2 y_4 y_5^2\\
  x_2 y_4^2 y_5 &=& x_2 y_4 y_5^2,
\end{eqnarray*}
which is a formidable system to digest by hand, especially when you
realize that it really does boil down to the three conceptual cases of
1 and 2 embedded at the same point (the other three placed freely
distance 1 around them), or 3,4,5 at the same point (with 1 and 2
free), or ***.  In any case, checking for duplicates, once a reduced
Gr\"obner basis is computed. is still a nontrivial task.

\section{Comments}
Chv\'atal et al. \cite[problem 21]{CKK72} posed the question in terms
of bounds on the number of vertices in a UDG, noting the lack of an
extant embedding for $H$ (in terms of projective planes). Hochberg
\cite{RH06} describes an algorithm for showing the impossibility of an
embedding for a given graph, unfortunately by experience not tractable
on $H$. Gerbracht \cite{EG06} found an analytic embedding for the
Harborth graph, the smallest known {\it non-crossing} UDG or
matchstick graph. He found a polynomial in one variable that
determines the finite set of possibilities for one point of the
embedding, from which all the rest are determined.

\bibliography{HeawoodUnitDistance}
\bibliographystyle{abbrv}

\end{document}